\documentclass[12pt]{article}
\usepackage{amsmath,amssymb,amscd,latexsym,euscript,mathrsfs}
\usepackage[matrix,arrow,curve]{xy}

\newcommand{\dom}{ {\rm dom}}
\newcommand{\cod}{ {\rm cod}}
\newcommand{\coLim}{\underrightarrow{\lim}}

\newcommand{\Lan}{{\rm Lan}}

\newcommand{\Set}{{\rm Set}}
\newcommand{\Ab}{{\rm Ab}}
\newcommand{\Ob}{{\rm Ob\,}}
\newcommand{\Mor}{{\rm Mor\,}}

\newcommand{\Tran}{{\rm Tran\,}}

\newcommand{\Imm}{{\rm Im\,}}
\newcommand{\Ker}{{\rm Ker\,}}

\newcommand{\NN}{{\,\mathbb N}}

\newcommand{\pt}{{\,\rm pt}}

\newcommand{\mC}{{\mathscr C}}
\newcommand{\mD}{{\mathscr D}}

\newcommand{\ZZ}{{\,\mathbb Z}}
\newcommand{\II}{{\,\mathbb I}}
\newcommand{\mA}{{\mathcal A}}

\newcommand{\mM}{{ M}}

\newcommand{\fS}{{\mathfrak{S}}}
\newcommand{\fF}{{\mathfrak{F}}}
\newcommand{\eproof}{{\Box}}

\newtheorem{theorem}{\bf Theorem}[section]
\newtheorem{lemma}[theorem]{\bf Lemma}
\newtheorem{proposition}[theorem]{\bf Proposition}
\newtheorem{corollary}[theorem]{\bf Corollary}
\newtheorem{definition}{\sc Definition}[section]
\newtheorem{example}[definition]{\sc Example}

\newtheorem{conjecture}{\sc Conjecture}

\def\leq{\leqslant}
\def\geq{\geqslant}

\begin{document}

\begin{center}
 {\large CUBICAL HOMOLOGY OF ASYNCHRONOUS \\ 
	TRANSITION SYSTEMS}
\\
A. A. Khusainov
\end{center}

\begin{abstract}
We show that a set with an action of a locally finite-dimensional free 
partially commutative monoid and the corresponding semicubical set have 
isomorpic homology groups. We build a complex of finite length for the 
computing homology groups of any asynchronous transition system with finite 
maximal number of mutually independent events.  
We give examples of computing the homology groups.
\end{abstract}

2000 Mathematics Subject Classification 18G10, 18G35, 55U10, 68Q10, 68Q85

Keywords: semicubical set, 
homology of small categories, free partially commutative monoid,  
trace monoid, asynchronous transition system.

\section*{Introduction}

By \cite{nie1996}, an {\em asynchronous transition system} 
$(S, s_0, E, I, \Tran)$ consists of arbitrary sets 
$E$  and $S$ with a distinguished element $s_0\in S$, 
 an irreflexive symmetric relation $I\subseteq E\times E$, and a subset
$\Tran\subseteq S\times E\times S$ satisfying the following axioms. 

\begin{enumerate}
\item for every $e\in E$, there are $s\in S$ and  $s'\in S$ such those $(s,e,s')\in \Tran$;
\item if $(s,e,s')\in \Tran$ and $(s,e,s'')\in \Tran$, then $s'=s''$;
\item for any pair $(e_1,e_2)\in I$ and triples $(s,e_1,s_1)\in \Tran$,
$(s_1,e_2,u)\in \Tran$ there exists $s_2\in S$ such that
$(s,e_2,s_2)\in \Tran$ and $(s_2,e_1,u)\in \Tran$.
\end{enumerate}
Elements $s\in S$ are called {\em states}, $e\in E$ {\em events}, 
$s_0$ is an {\em initial state} and $I\subseteq E\times E$  
{\em independence relation}. Triples $(s,e,s')\in \Tran$ are {\em transitions}.

Asynchronous transition systems are introduced by
Mike Shields \cite{shi1985} and Marek Bednarczyk \cite{bed1988} for 
a simulation of concurrent computational systems. 
The application of the partially commuting 
in parallel programming belongs to Antoni Mazurkievicz.
In  \cite{X20032}, it was proposed to consider an asynchronous transition system 
as a pointed set with action of a free partially commutative monoid.
It allows to introduce \cite{X20032} and find approach to studying homology 
groups of asynchronous transition systems \cite{X20042}.
Erik Goubault \cite{gou1995} and Philippe Gaucher \cite{gau2000}
have defined homology groups for higher dimensional automata, 
which are other models of parallel computational systems.  
The our main result  (Theorem \ref{mainres}) 
shows that in some finiteness conditions, 
the homology groups of asynchronous transition systems 
and of corresponding semiregular higher dimensional automata 
are isomorpic
(see Corollary \ref{corres}).
In \cite{X20042}, it was built an algorithm 
to computing the first integer homology group 
of any asynchronous transition system. 
In \cite[òåîðåìà 1]{X20082}, by a resolution of Ludmila Polyakova
\cite{pol2007}, it was built a complex for the computations 
of homology groups in the case of finite asynchronous transition system.
We will build a complex for the computing 
the homology groups of an asynchronous transition system
  without infinite sets of mutually 
independent events (Corollary \ref{main3}).
If the maximal number of mutually independent events is finite, then 
this complex has a finite length.

\section{Homology of categories and semicubical sets}

Let $\mA$ be  a category. 
Denote by $\mA^{op}$ the opposite category.
Given $a,b\in \Ob\mA$, let   
$\mA(a,b)$ be the set of all morphisms $a\rightarrow b$.  
For a small category  $\mC$, denote by $\mA^{\mC}$  the category 
of functors $\mC\rightarrow \mA$ and natural transformations.

Throughout this paper let 
$\Set$ be the category of sets and maps,  
$\Ab$ the category of abelian groups and homomorphisms,
$\ZZ$ the set or additive group of integers,
$\NN$ the set of nonnegative integers or 
 free monoid  $\{1, a, a^2, \cdots \}$ generated by one element.
For any family of abelian groups $\{A_j\}_{j\in J}$, the direct sum 
is denoted by $\bigoplus\limits_{j\in J}A_j$. Elements of direct summands 
is written as pairs $(j,g)$ with $j\in J$ and $g\in A_j$.
If $A_j= A$ for all $j\in J$, then the direct sum is denoted by
 $A^{(J)}$.

\subsection{Semicubical sets}

Suppose that $\II=\{0,1\}$ is the set ordered by $0<1$. 
For integer $n\geq 0$, let $\II^n$ be the Cartesian power 
of $\II$. Denote by
 $\Box_+$ the category of partially ordered sets 
$\II^n$ 
and maps, which can be decomposed into compositions of the increasing maps 
$\delta_i^{k,\varepsilon}: \II^{k-1}\rightarrow \II^k$, 
$1 \leq i \leq k$, $\varepsilon\in \II$ defined by
 $\delta_i^{k,\varepsilon}(x_1, \cdots, x_{k-1})=
(x_1, \cdots, x_{i-1}, \varepsilon, x_i, \cdots, x_{k-1}).$

A {\em semicubical set} \cite{X20081} is any functor 
$X: \Box_{+}^{op} {\rightarrow}\Set$. Morphisms are defined as 
natural transformations. 
Since every morphism $f: \II^m\rightarrow \II^n$ of the category $\Box_+$ 
has the canonical decomposition
$f=\delta_{j_{n-m}}^{n,\varepsilon_{n-m}}\cdots\delta_{j_1}^{m+1,\varepsilon_1}$ 
such that
 $1\leq j_1 < \cdots < j_{n-m}\leq n$, a functor $X$ is defined by 
values  
 $X_n=X(\II^n)$ on objects and  
$\partial_i^{k,\varepsilon}=X(\delta_i^{k,\varepsilon})$ on morphisms.
Hence a semicubical set may be given 
as a pair
$(X_n, \partial_i^{n,\varepsilon})$ consisting of a sequnce of sets
$(X_n)_{n\in \NN}$ and a family of maps 
$\partial_i^{n,\varepsilon}: X_n \rightarrow X_{n-1}$ defined  
for $1\leq i\leq n$, $\varepsilon\in \{0,1\}$, and satisfying the condition
$$
\partial_i^{n-1,\alpha}\circ \partial_j^{n,\beta} =
\partial_{j-1}^{n-1,\beta}\circ \partial_i^{n,\alpha}~,
\mbox{ for } \alpha,\beta \in \{0,1\}, n\geq 2, \mbox{ and } 1\leq i< j\leq n.
$$

For example, any directed graph given by a pair of maps 
$\dom, \cod: X_1\rightarrow X_0$ assigning to every arrow 
its source and target, we can considered as the semicubical set 
with $\partial_1^{1,0}=\dom$, 
$\partial_1^{1,1}= \cod$, and $X_n=\emptyset$ for all $n\geq 2$.
For $n\geq 2$, the maps $\partial_i^{n,\varepsilon}$ are empty.
By \cite{X20081}, {\em cubical sets} \cite{kac2004} provide 
examples of semicubical sets.
Similarly, we can define  {\em semicubical objects} 
$(X_n, \partial_i^{n,\varepsilon})$ 
in an arbitrary category  $\mA$.

\subsection{Homology of small categories}

Homology groups of small categories 
with coefficients in functors into the category of Abelian groups will be considered in this 
subsection.

\smallskip
\noindent
{\bf Homology of categories and derived functors of the colimit.} 

\begin{definition}\label{defhomolcat}
Let $\mC$ be a small category and 
 $F:\mC\rightarrow \Ab$ a functor. Denote by 
$C_*(\mC, F)$ a chain complex of Abelian groups
\begin{displaymath}
 C_n({\mC},F) = \bigoplus_{c_0 \rightarrow \cdots \rightarrow c_n}
F(c_0), \quad n \geq 0,
\end{displaymath}
with differentials
$d_n= \sum\limits_{i=0}^{n}(-1)^i d^n_i: 
C_n(\mC,F) \rightarrow C_{n-1}(\mC,F)$ defined for $n>0$ by the {\em face operators}
\begin{multline*}
d^n_i(c_0 \stackrel{\alpha_1}\rightarrow c_1 \stackrel{\alpha_2}\rightarrow
\cdots \stackrel{\alpha_{n}}\rightarrow c_{n}, a)=\\
\left\{
\begin{array}{ll}
(c_1 \stackrel{\alpha_2}\rightarrow \cdots \stackrel{\alpha_n}\rightarrow c_n, 
F(c_0\stackrel{\alpha_1}\rightarrow c_1)(a) ) ~~, 
& \mbox{if}~ i = 0\\
(c_0 \stackrel{\alpha_1}\rightarrow \cdots \stackrel{\alpha_{i-1}}\rightarrow c_{i-1} 
\stackrel{\alpha_{i+1}\alpha_{i}}\rightarrow c_{i+1} \stackrel{\alpha_{i+2}}\rightarrow 
\cdots \stackrel{\alpha_n}\rightarrow c_n, a) \quad , & \mbox{if} ~ 1 \leq i \leq n-1\\
(c_0 \stackrel{\alpha_1}\rightarrow \cdots \stackrel{\alpha_{n-1}}\rightarrow c_{n-1}, 
a ) ~~, & \mbox{if}~ i=n
\end{array}
\right.
\end{multline*}
The quotient groups $H_n(C_*(\mC,F))=\Ker(d_n)/\Imm(d_{n+1})$ are called the
{\em $n$-th homology groups of the category  $\mC$ with coefficients in $F$}.
\end{definition}
It is known \cite{X20081} that the functors  
$H_n(C_*(\mC,-)): \Ab^\mC \rightarrow \Ab$ are natural isomorphic 
to left derived of the colimit functor
 $\coLim^\mC: \Ab^\mC \rightarrow \Ab$. 
So we will be denote them by $\coLim_n^\mC$.

For any small category $\mC$,
denote by 
 $\Delta_\mC\ZZ$, or $\Delta\ZZ$ shortly, the functor $\mC\rightarrow \Ab$
with constant values  
 $\ZZ$ on objects and $1_\ZZ$ on morphisms.
The values of left satellites $\coLim_n^\mC\Delta_\mC\ZZ$ of the colimit functor on 
 $\Delta_\mC\ZZ$ are called the {homology groups of the category $\mC$} and denoted by $H_n(\mC)$.
It follows from Eilenberg's Theorem \cite[Appl. 2]{gab1967} that 
 homology groups of the geometric realization of the nerve of $\mC$
are isomorphic to $H_n(\mC)$.
For example, if the category denoted by $\pt$ consists of unique object and the identity morphism, 
then $H_n(\pt)=0$ for all $n>0$ and $H_0(\pt)=\ZZ$.

\noindent
{\bf Coinitial functors.}  
Let $\mC$ be a small category.
If $H_n(\mC)=0$ for $n>0$ and $H_0(\mC)\cong \ZZ$, then $\mC$ is called {\em acyclic}.
Let  $S: \mC \rightarrow \mD$ be a functor into an arbitrary category $\mD$.
For any $d\in \Ob(\mD)$,  a {\em fibre} 
(or {\em comma-category} \cite{mac1998}), denoted by
$S/d$, is the category whose objects $\Ob(S/d)$ consists of pairs $(c,\alpha)$ with
$c\in \Ob(\mC)$, and $\alpha\in \mD(S(c),d)$.
Morphisms in $S/d$ are triples $(f, \alpha_1, \alpha_2)$ of morphisms
$f\in \mC(c_1,c_2)$, $\alpha_1\in \mD(S(c_1),d)$, and $\alpha_2\in \mD(S(c_2),d)$ 
satisfying $\alpha_2\circ S(f)=\alpha_1$.
A {\em forgetful functor $Q_d: S/d\rightarrow \mC$ of the fibre} is defined by
$Q_d(c,\alpha)=c$ on objects and $Q_d(f, \alpha_1, \alpha_2)=f$
on morphisms.
If the functor $S$ is a full inclusion $\mC\subseteq \mD$, then $S/d$
is denoted by $\mC/d$.
\begin{definition}
 If the category $S/d$ is acyclic for every object $d\in \mD$,
then à functor $S: \mC \rightarrow \mD$ 
is called {\em strongly coinitial}.
\end{definition}

By
\cite[Theorem 2.3]{obe1968}, it can be proved, that a functor
$S:\mC\rightarrow\mD$ between small categories is strongly 
coinitial if and if the canonical morphisms 
$\coLim_n^{\mC^{op}}(F\circ S^{op})\rightarrow \coLim_n^{\mD^{op}}F$
are isomorhisms for all functors and $n\geq 0$ \cite[Proposition 1.4]{X20083}.

\subsection{Homology of semicubical sets}

Let  $X\in \Set^{\Box_+^{op}}$ be a semicubical set. 
The {\em category of singular cubes}
$h_*/X$ is the fibre of the Yoneda Embedding $h_*: \Box_+\rightarrow \Set^{\Box_+^{op}}$ 
over $X$. 
Consider a category  
$\Box_+/X$ whose objects are elements $\sigma\in \coprod\limits_{n\in\NN}X_n$. 
Morhisms from 
$\sigma\in X_m$ to $\tau\in X_n$ are given by triples 
$(\alpha, \sigma, \tau)$,  $\alpha\in \Box_+(\II^m,\II^n)$, satifying 
 $X(\alpha)(\tau)=\sigma$. The categories  $h_*/X$ and $\Box_+/X$ 
are isomorphic and we will identify their.
A {\em homological system on a semicubical set
$X$} is any functor $F:(\Box_+/X)^{op}\rightarrow \Ab$.

\medskip
\noindent
{\bf Homology groups with coefficients in a homological system.} 
We turn to study homology groups  $\coLim_n^{(\Box_+/X)^{op}}F$ 
of the category of singular cubes with coefficients in a 
homological system $F$ on $X$. 
Consider abelian groups
$C_n(X,F)=\bigoplus\limits_{\sigma\in X_n}F(\sigma)$. Define  
{\em boundary operators} $d_i^{n,\varepsilon}: C_n(X,F)\rightarrow C_{n-1}(X,F)$ as 
homomorphisms such those the following diagrams commute 
for $1\leq i\leq n$ and 
$\varepsilon \in \II= \{0, 1\}$,
$$
\begin{CD}
\bigoplus\limits_{\sigma\in X_n}F(\sigma) @>d_i^{n,\varepsilon}>> 
		\bigoplus\limits_{\sigma\in X_{n-1}}F(\sigma)\\
@A in_\sigma AA @A in_{X(\delta^{n,\varepsilon}_i)(\sigma)}AA\\
F(\sigma) 
@>F(\delta_i^{n,\varepsilon}, X(\delta^{n,\varepsilon}_i)(\sigma), \sigma)>> 
F(X(\delta^{n,\varepsilon}_i)(\sigma))
\end{CD}
$$

\begin{definition}
Let $F: (\Box_+/X)^{op}\rightarrow \Ab$
be a homological system on a semicubical set $X$. 
{\em Homology groups $H_n(X,F)$ with coefficients in $F$}
are $n$-th homology groups of the complex $C_*(X,F)$ consisting of the groups 
$C_n(X,F)= \bigoplus\limits_{\sigma\in X_n}F(\sigma)$ 
and differentials $d_n=\sum\limits_{i=1}^n (-1)^i (d^{n,1}_i-d^{n,0}_i)$. 
\end{definition}
Groups $H_n(X,\Delta\ZZ)$ are called $n$-th {\em integer}
 homology groups.

By \cite[òåîðåìà 4.3]{X20081}, 
for any semicubical set $X$ and a homological system 
$F$ on $X$, there is isomorphisms 
$\coLim_n^{(\Box_+/X)^{op}}F  \cong  H_n(X,F)$ for all $n\geq 0$.
It follows that homology groups of cubical sets studied in  \cite{kac2004} and  \cite{kac2003} 
are isomorphic to homology groups of the corresponding 
semicubical sets with constant homological systems.

\section{Homology of free partially commutative\\
 monoids}

In this section, we study a factorization category 
and semicubical set concerned with a free partially commutative monoid.
We will prove that homology groups 
of the free partially commutative monoid 
and the corresponding semicubical set are isomorphic. 
It allows us to build a complex to the computing 
the homology groups of a free partially commutative monoid.

Each monoid $M$ will be considered as a category with 
 $\Ob M= \{M\}$ consisting of the unique object and $\Mor M= M$.
It has an influence on denotations and terminology. 
In particular, a right $M$--set $X$, with the action 
$(x,\mu)\mapsto x\cdot\mu$ for $x\in X$ and $\mu\in M$ is considered and denoted 
as the functor ${X}: \mM^{op}\rightarrow \Set$ defined by ${X}(\mu)(x)= x\cdot\mu$. 
Morphisms of right $M$-sets are natural transformations.

\subsection{Homology of a free finitely generated commutative 
monoid}

{\bf Category of factorizations.} 
 Suppose that $\mC$ is a small category. 
If a morphism $f\in \Mor\mC$ belongs to $\mC(a,b)$, then we write $\dom f = a$
and $\cod f = b$. A {\em category of factorizations} $\fF\mC$ 
\cite{bau1985} 
has the set of objects $\Ob(\fF\mC)= \Mor(\mC)$ and the sets of morphisms 
$\fF\mC(\alpha,\beta)$ 
consisting of all pairs $(f,g)$ of $f\in \mC(\dom\beta, \dom\alpha)$ and 
$g\in \mC(\cod\alpha, \cod\beta)$ satisfying 
 $g\circ\alpha\circ f = \beta$. 

The composition of morphisms $\alpha\stackrel{(f_1,g_1)}\rightarrow \beta$ and 
$\beta\stackrel{(f_2,g_2)}\rightarrow \gamma$
 is defined by 
$\alpha\stackrel{(f_1\circ f_2,g_2\circ g_1)}\longrightarrow \gamma$.
The identity of an object $a\stackrel\alpha\rightarrow b$ of the category  
$\fF\mC$ consists of the pair of identity morphisms 
$ \alpha \stackrel{(1_a, 1_b)}\longrightarrow \alpha$.
We will study the category of factorizations of a monoid 
considered as a category with an unique object.
 Denote by
$\fF\mC \stackrel{\cod}\rightarrow \mC$ the functor that assigns to each  
$a\rightarrow b$ its the codomain $b$ and to each morphism   
$\alpha \stackrel{(f,g)}\rightarrow \beta$ the morphism   
$g: \cod(\alpha)\rightarrow \cod(\beta)$. By \cite[Lemma 1.9]{X20083}, this functor is strongly 
coinitial.

\begin{lemma}\label{strcof}
Let $\NN=\{1, a, a^2, \ldots\}$ be the free monoid generated by 
 $a$ and let $T=\{1,a\}$. Suppose that $\fF{T}\subset \fF\NN$ is  
full subcategory with the set of objects $T$. Then the inclusion $\fF{T}\subset \fF\NN$ is 
strongly coinitial.
\end{lemma}
{\sc Proof.} 
The category $\fF{T}$ consists of two objects and two morphisms
$1\stackrel{(1,a)}\rightarrow a$, $1\stackrel{(a,1)}\rightarrow a$
except the identity morphisms. It easy to see that for any integer  
 $p\geq 0$, 
the comma-category $\fF{T}/a^p$ is the poset
\begin{multline*}
(1,1,a^p) < (1,a,a^{p-1}) >  (a^1,1, a^{p-1}) < 
\cdots \\ 
 \cdots > (a^s,1, a^t) < (a^s, a, a^{t-1}) > (a^{s+1}, 1, a^{t-1}) 
< \cdots > (a^p,1,1)
\end{multline*}
Since the geometric realization of its nerve 
is homeomorphic to the unit segment,
$H_q(\fF{T}/a^p)\cong H_q(pt)$ for all $q\geq 0$. 
\hfill $\eproof$

Suppose that $\NN$ is the free monoid generated by $a$. For $n\geq 1$, denote 
$a_1=(a, 1, \cdots, 1)$, $a_2=(1, a, 1, \cdots, 1)$, $\cdots$, 
$a_n=( 1, \cdots, 1, a)$. Consider the subset $T^n\subset \NN^n$ consisting of 
finite products $a_{i_1}a_{i_2}\cdots a_{i_k}$ where 
$1\leq i_1< i_2 < \cdots < i_k\leq n$. 
Let $\fF T^n$ be the full subcategory of  $\fF\NN^n$  with the class of objects $T^n$. 
For every $\alpha=(a^{p_1}, a^{p_2}, \cdots, a^{p_n})$, the comma-category 
 $\fF T^n/\alpha$ is isomorphic to 
$(\fF T/a^{p_1})\times \cdots
 \times(\fF T/a^{p_n})$. It follows from the K\"unneth formula \cite[Lemma 1.16]{X20083}
and Lemma \ref{strcof} the following.

\begin{lemma}\label{confn}
The inclusion $\fF T^n \subset \fF\NN^n$ is strongly coinitial.
\end{lemma}

\medskip
\noindent
{\bf A semicubical set of a free finitely generated 
commutative monoid.}
Let  
$a_1$, $a_2$, $\dots$, 
$a_n$ the above generators of $\NN^n$.
 Consider the semicubical set  $T^n_*$ consisting of the subsets
$$
T^n_k = \{a_{i_1}a_{i_2}\cdots a_{i_k} : 
1\leq i_1 < i_2 < \cdots < i_k \leq n \}
$$
and maps $T^n_{k-1} ~~ { {{\partial^{k,0}_s}\atop\longleftarrow}\atop
{\longleftarrow\atop{\partial^{k,1}_s}} } ~~ T^n_k ~$, $1\leq s\leq k$ defined as follows.
$$
\partial_s^{k,0}(a_{i_1}\cdots a_{i_k}) = \partial_s^{k,1}(a_{i_1}\cdots a_{i_k}) =
a_{i_1}\cdots a_{i_{s-1}}\widehat{a_{i_s}}a_{i_{s+1}} \cdots a_{i_k}
$$
Here  $a_{i_1}\cdots a_{i_{s-1}}\widehat{a_{i_s}}a_{i_{s+1}} \cdots a_{i_k}$ is the word 
obtained by removing the symbol $a_{i_s}$.
Objects of the category 
$\Box_+/T^n_*$ may be considered as pairs $(k, \sigma)$ where $\sigma\in T^n_k$. Every $\sigma\in T^n_k$
 has an unique decomposition  
$a_{i_1}\cdots a_{i_k}$ such that 
$1\leq i_1 < \cdots < i_k\leq n$.  And so the objects $(k,\sigma)$ may 
be identify with elements
 $a_{i_1}\cdots a_{i_k}\in T^n$. 
Morphisms in $\Box_+/T^n_*$ are triples  
$$
(\delta: \II^m \rightarrow \II^k, a_{j_1}\cdots a_{j_m}, 
a_{i_1}\cdots a_{i_k}),$$
such those $T^n(\delta)(a_{j_1}\cdots a_{j_m})= a_{i_1}\cdots a_{i_k}$.

We will establish a functor $\fS: \Box_+/T^n_* \rightarrow  \fF T^n$.
Toward this end, define 
$\fS(a_{i_1}\cdots a_{i_k})=a_{i_1}\cdots a_{i_k}$ on objects.
Every morphism of the category $\Box_+/T^n_*$ has a decomposition 
$(\delta_s^{k,\varepsilon}, a_{i_1}\cdots \widehat{a_{i_s}} \cdots a_{i_k}, 
a_{i_1}\cdots a_{i_k})$, $\varepsilon\in \{0,1\}$.
Hence, it is enough to define the values 
\begin{multline}\label{cubtofrac}
 \fS(\delta_s^{k,0}, a_{i_1}\cdots \widehat{a_{i_s}} \cdots a_{i_k}, 
a_{i_1}\cdots a_{i_k}) = 
(a_{i_s},1): a_{i_1}\cdots \widehat{a_{i_s}} \cdots a_{i_k}
\rightarrow 
a_{i_1}\cdots a_{i_k}\\ 
\fS(\delta_s^{k,1}, a_{i_1}\cdots \widehat{a_{i_s}} \cdots a_{i_k}, 
a_{i_1}\cdots a_{i_k}) = 
(1, a_{i_s}): a_{i_1}\cdots \widehat{a_{i_s}} \cdots a_{i_k}
\rightarrow 
a_{i_1}\cdots a_{i_k}
\end{multline} 
It is easy to see that the map $\fS$
has the unique functorial extension. 

For each $\sigma= a_{i_1}\cdots a_{i_k}$,
 the category $\fS/\sigma$ has the terminal object. Therefore, 
 the following assertion holds.

\begin{lemma}\label{miscpremain}
The functor $\fS: \Box_+/T^n_* \rightarrow  \fF T^n$ is strongly coinitial.
\end{lemma}

\subsection{Homology of a free partially commutative monoid 
with coefficient in a right module}

\begin{definition}
Let $E$ be a set. Suppose that $I\subseteq E\times E$ is an irreflexive 
symmetric relation on $E$. 
Monoid
 given by the set of generators $E$ and the relations $ab=ba$ for all $(a,b)\in I$ 
is called {\em free partially commutative} and denoted by  $M(E,I)$.
If $(a,b)\in I$, then the elements $a,b \in E$ are called the
{\em commuting generators}.
\end{definition}
Our definition is more general than it is given in \cite{die1997}. 
We do not demand that the set $E$ is finite.

For any graph, a subgraph is called a {\em $n$-clique} if that is 
isomorphic to the complete graph $K_n$. 
A {\em clique} is a subgraph which is equal to a $n$-clique for some 
cardinal number $n\geq 1$.
Let $M(E,I)$ be a free partially commutative monoid
given by a set of generators $E$ and relations $ab=ba$ for all $(a,b)\in I$. 
Denote by $V$ the set of all maximal cliques of its indepedence graph.
For every $v\in V$, denote by
 $E_v\subseteq E$ the set of  vertices of $v$. 
The set $E_v$ is a maximal subset of $E$
consisting of mutually commuting elements.
The set  $E_v$ generate the maximal commutative submonoid 
$M(E_v)\subseteq M(E,I)$.
The monoid $M(E,I)$ is called {\em locally finite-dimensional} if the sets 
$E_v$ are finite for all $v\in V$.
This property holds if and only if 
the independence graph does not contain infinite cliques.

\medskip
\noindent
{\bf The coinitial subcategory of a category of factorizations.}
As above, V is the set of maximal cliques in the independence graph of $M(E,I)$.
\begin{proposition}
Suppose for $v\in V$ that $T_v\subset M(E_v)$  is the subset 
of products
$a_1 a_2 \cdots a_n$ of distinct elements  
$a_j\in E_v$, $1\leq j\leq n$. Here $n$ may be taken finite values  $\leq \vert E_v \vert$.  
It is supposed that the product equals  $1\in T_v$ for $n=0$.
If $M(E,I)$ is locally finite-dimensional, then the inclusion
 $\bigcup\limits_{v\in V} \fF T_v \subset \fF M(E,I)$ is strongly coinitial.
\end{proposition}
{\sc Proof.} The composition of strongly coinitial functors is strongly coinitial.
 Since the inclusion  
$\bigcup\limits_{v\in V} \fF M(E_v) \subseteq \fF M(E,I)$ is strongly 
coinitial by \cite[Theorem 2.3]{X20083}, it is enough to show that the inclusion 
$\bigcup\limits_{v\in V} \fF T_v \subset\bigcup\limits_{v\in V} \fF M(E_v)$
is strongly coinitial.
For each $\alpha \in \bigcup\limits_{v\in V} \fF M(E_v)$, there is  
 $w\in V$ such that $\alpha \in \fF M(E_w)$. All divisors of $\alpha$
belong to $M(E_w)$. Hence
$\bigcup\limits_{v\in V} \fF T_v/\alpha = \fF T_w/\alpha$.
By Lemma \ref{confn}, the inclusion  
$\fF T_w\subset \fF M(E_w)$ is strongly coinitial. Therefore 
$H_q(\bigcup\limits_{v\in V} \fF T_v/\alpha)\cong H_q(pt)$.
\hfill $\eproof$

\medskip
\noindent
{\bf Cubical homology of free partially commutative monoids.}
For an arbitrary set $E$ with an irreflexive symmetric relation 
 $I\subseteq E\times E$, we construct a semicubical set
$T(E,I)$ depending on a some total ordering relationship $\leq$ on $E$.
Toward this end for any integer $n>0$, we define 
 $T_n(E,I)$ as the set of all tuples $(a_1, \cdots, a_n)$ 
of mutually commuting elements $a_1<\cdots < a_n$ in $E$,
$$
T_n(E,I) = \{(a_1, \cdots, a_n): (a_1 < \cdots < a_n) \& 
( 1\leq i<j\leq n \Rightarrow (a_i, a_j)\in I) \}.
$$
The set $T_0(E,I)$ consists of unique empty word $1$.
Maps $\partial^{n,\varepsilon}_i: T_n(E,I)\rightarrow T(E,I)_{n-1}$
for $1\leq i\leq n$ act as
\begin{equation}\label{cubeq}
\partial^{n,0}_i(a_1, \cdots,  a_n) = \partial^{n,1}_i(a_1, \cdots, a_n) = 
(a_1, \cdots,  \widehat{a_i},  \cdots, a_n)
\end{equation}

It easy to see that $T(E,I)$ is equal to union of the semicubical sets 
$(T_v)_*$ defined by 
$$
(T_v)_n = \{(a_1, \cdots, a_n) \in T_n(E,I)~:~ ( 1\leq i\leq n \Rightarrow
a_i\in E_v)\},
$$
where $E_v\subseteq E$ are the maximal subsets of mutually commuting 
generators of $M(E,I)$.
Face operators 
$(T_v)_n \stackrel{\partial^{n,\varepsilon}_i}\rightarrow (T_v)_{n-1}$
act for $1\leq i \leq n$ and $\varepsilon \in \{0,1\}$
 by (\ref{cubeq}).

Let $\fS: \Box_+/\bigcup\limits_{v\in V} (T_v)_* \rightarrow 
\bigcup\limits_{v\in V} \fF T_v$ be the functor assigning to every singular 
qube  $(a_1,\cdots, a_n) \in \bigcup\limits_{v\in V} (T_v)_n$
the object $a_1\cdots a_n$. The functor $\fS$ acts on morphisms 
by the equation (\ref{cubtofrac}).

For each functor
$F: M(E,I)^{op} \rightarrow \Ab$,
denote by $\overline{F}$
a homological system on $T(E,I)$ defined as the composition 
$$
(\Box_+/T(E,I))^{op} \stackrel{\fS^{op}}\longrightarrow 
{\bigcup\limits_{v\in V}(\fF T_v)^{op}} \subset (\fF M(E,I))^{op}
\stackrel{\cod^{op}}\longrightarrow M(E,I)^{op}
\stackrel{F}\rightarrow  \Ab
$$

We will suppose below that $E$ is a totally ordered set.
\begin{proposition}\label{leechiscubic}
Let $M(E,I)$ be a locally finite-dimensional free partially commutative monoid and  
$F: M(E,I)^{op} \rightarrow \Ab$ a functor.
The homology groups of $M(E,I)$ are isomorphic to the cubical homology groups 
$$
\coLim_n^{M(E,I)^{op}}F \cong H_n(T(E,I),\overline{F}), \quad n\geq 0.
$$
\end{proposition}
{\sc Proof.}  
The image of $\fS|_{\Box_+/(T_v)_*}$ 
is contained in the category $\fF T_v$ and defines a functor 
which we denote by
$\fS_v: \Box_+/ (T_v)_* \rightarrow \fF T_v$.
For any 
$\alpha \in \bigcup\limits_{v\in V} \fF T_v$, there exists $w\in V$ such that
$\alpha \in \fF T_w$. It follows that there is an isomorphism
$\fS/\alpha \cong \fS_w/\alpha$.
By Lemma \ref{miscpremain}, for every free finitely generated 
monoid $M(E_v)$, the functor $\fS_v: \Box_+/ (T_v)_* \rightarrow  \fF T_v$ 
is strongly coinitial. Consequently, the category  
$\fS_v/\alpha$  is acyclic. It follows that  
$\fS/\alpha$  is acyclic. Thus $\fS$ is strongly coinitial.
Using the strong coinitiality of the inclusion 
$\bigcup\limits_{v\in V} \fF T_v \subseteq \fF M(E,I)$, we have 
the isomorphisms  
$$
\coLim_n^{M(E,I)^{op}} F \cong 
\coLim_n^{(\Box_+/\bigcup\limits_{v\in V} 
(T_v)_*)^{op}}~\overline{F}
$$
\hfill $\eproof$

It allows us to construct the following complex for the computing the groups 
$\coLim_n^{M(E,I)^{op}}F$.

\begin{corollary}\label{hrm}
Suppose that $M(E,I)$ is a locally finite-dimensional free partially commutative monoid. 
Then for any right $M(E,I)$-module $G$, the groups $H_n(M(E,I)^{op},G)$
are isomorphic to homology groups of the complex 
\begin{multline*}
0 \leftarrow G \stackrel{d_1}\leftarrow \bigoplus\limits_{a_1\in T_1(E,I)} G
\stackrel{d_2}\leftarrow \bigoplus\limits_{(a_1,\ldots, a_n)\in T_n(E,I)} G
\leftarrow  \cdots \\
\cdots \leftarrow  
\bigoplus\limits_{(a_1,\ldots, a_{n-1})\in T_{n-1}(E,I)} 
G
\stackrel{d_n}\longleftarrow 
\bigoplus\limits_{(a_1,\ldots, a_n)\in T_n(E,I)} 
G \leftarrow \cdots~,
\end{multline*}
the $n$-th member of that for each $n\geq 0$ equals a direct sum of backup copies of the abelian group 
$G$ 
in all $n$-tuples of mutually commuting elements $a_1<a_2<\cdots<a_n$ in $E$
where the differentials are defined by 
\begin{equation}\label{diff}
d_n(a_1, \cdots, a_n, g) = 
\sum\limits_{s=1}^n (-1)^s (a_1, \cdots, \widehat{a_s} , \cdots, a_n,
G(a_s)(g)-g)
\end{equation}
\end{corollary}

\section{Homology of sets with an action of a free partially 
commutative monoid}

This section is devoted to homology groups of right $M(E,I)$-sets $X$
with coefficients in functors $F: (M(E,I)/X)^{op}\rightarrow \Ab$.
We show that these groups may be studied 
as the homology groups of the monoid $M(E,I)$.
We prove the main result of this paper about the isomorphism of homology groups 
of any  $M(E,I)$-set $X$ and the corresponding semicubical set $Q_*X$. 
This result is applied for the computing 
homology groups of state spaces in the simplest cases.

\subsection{Homology of right $M(E,I)$-sets} 

Let  $M$ be a monoid and  $X\in \Set^{M^{op}}$ a right $M$-set. 
Suppose that $F: (M/X)^{op}\rightarrow \Ab$ is a functor. 
Consider a functor $S=Q^{op}: (M/X)^{op}\rightarrow M^{op}$ 
opposite to the forgetful functor of the fibre. 
Let $\Lan^S F: M^{op}\rightarrow Ab$ be the left Kan extension \cite{mac1998}
of the functor $F$ along to $S$.
Objects of the category $(M/X)^{op}$ may be considered as elements
 $x\in X$. Morphisms between them are triples $x\stackrel{\mu}\rightarrow y$ such those
 $\mu\in M$ and $x\cdot\mu=y$.
We get the following assertion 
by replacing  a monoid instead of the category 
in \cite[Proposition 3.7]{X20081}.

\begin{lemma}\label{misc2}
A right $M$-module $\Lan^S F$ is the Abelian group
$\bigoplus\limits_{x\in X}F(x)$ with the action on $(x,f)$ with $x\in X$ and 
$f\in F(x)$
 defined by
$$
(x,f)\mu= \Lan^S F(\mu)(x, f)= (x\cdot\mu, F(x\stackrel{\mu}\rightarrow  x\cdot\mu)(f)).
$$ 
There are isomorphisms
$\coLim_n^{(M/X)^{op}}F\cong \coLim_n^{M^{op}}\Lan^S{F}$ for all $n\geq 0$.
\end{lemma}

Let $X$ be a right $M(E,I)$-set.
Consider an arbitrary total ordering on $E$.
Define sets  
$$
Q_n X = \{(x, a_1, \cdots, a_n): a_1<\cdots< a_n \&
(1\leq i < j \leq n \Rightarrow (a_i,a_j)\in I) \}.
$$
In particular $Q_0 X= X$. Define the maps 
$$
 Q_n X {{{\partial^{n,0}_i}\atop{\longrightarrow}}\atop
	{{\longrightarrow}\atop{\partial^{n,1}_i}}}  Q_{n-1}X~,  
\quad n\geq 1~, ~~1\leq i\leq n~, 
$$
by 
$$
{\partial^{n,\varepsilon}_i}(x,a_1, \cdots, a_n) = 
(x\cdot a_i^{\varepsilon}, a_1, \cdots, \widehat{a_i}, \cdots, a_n)~, \varepsilon\in \{0,1\},
$$ 
where $a_i^0=1$ è $a_i^1=a_i$.
\begin{lemma}
The sequence of the sets $Q_n X$ and the family of maps 
${\partial^{n,0}_i}$, ${\partial^{n,1}_i}$
make up a semicubical set.
\end{lemma}

For any functor $F: (M(E,I)/X)^{op}\rightarrow \Ab$,
we build the homological system 
$\overline{F}: (\Box_+/Q_*X)^{op}\rightarrow \Ab$ as follows. 
Define 
 $\overline{F}(x, a_1, \cdots, a_n)= F(x)$ on objects and 
$\overline{F}(\delta_i^{n,\varepsilon},  \partial_i^{n,\varepsilon}(\sigma), \sigma)
= F(x\stackrel{a_i^{\varepsilon}}\rightarrow xa_i^{\varepsilon})$
on morphisms for $\sigma=(x, a_1, \cdots, a_n)$.

\begin{theorem}\label{mainres}
Let $M(E,I)$ be a locally finite-dimensional free partially 
commutative monoid and 
$X$ a right $M(E,I)$-set.
Suppose $F: (M(E,I)/X)^{op} \rightarrow \Ab$ is a functor and
$\overline{F}: (\Box_+/Q_*X)^{op}\rightarrow \Ab$ the corresponding homological 
system of Abelian groups.
Then $\coLim_n^{(M(E,I)/X)^{op}}F\cong H_n(Q_*X,\overline{F})$ 
for all  $n\geq 0$.
In other words, $\coLim_n^{(M(E,I)/X)^{op}}F$ are isomorphic to homology groups 
of the complex
\begin{multline*}
0 \leftarrow 
\bigoplus\limits_{x\in Q_0 X} F(x) \stackrel{d_1}\leftarrow 
\bigoplus\limits_{(x,a_1)\in Q_1 X} F(x)
\stackrel{d_2}\leftarrow \bigoplus\limits_{{(x, a_1, a_2)\in Q_2 X}
} F(x)
\leftarrow  \cdots \\
\cdots \leftarrow  
\bigoplus\limits_{(x, a_1, \cdots, a_{n-1})\in Q_{n-1}X 
 } 
F(x)
\stackrel{d_n}\longleftarrow 
\bigoplus\limits_{{(x, a_1, \cdots, a_n)\in Q_n X} } 
F(x) \leftarrow \cdots~,
\end{multline*}
where $d_n(x,a_1, \cdots, a_n,f) = $
\begin{multline*}
\sum_{s=1}^n(-1)^s (
(x\cdot a_s, a_1, \cdots, \widehat{a_s}, \cdots, a_n, 
F(x\stackrel{a_s}\rightarrow x\cdot a_s)(f))\\
- (x, a_1, \cdots, \widehat{a_s}, \cdots, a_n, f) )
\end{multline*}
\end{theorem}
{\sc Proof.}
Applying Lemma \ref{misc2} to the monoid $M=M(E,I)$ gives 
 the following complex for 
a computation of homology groups with coefficients 
in  $\Lan^S F$ by Corollary \ref{hrm}.

\begin{multline*}
0 \leftarrow \Lan^S F \stackrel{d_1}\leftarrow \bigoplus\limits_{a_1\in T_1(E,I)} \Lan^S F
\stackrel{d_2}\leftarrow \bigoplus\limits_{(a_1, a_2)\in T_2(E,I)} \Lan^S F
\leftarrow  \cdots 
\\
\cdots \leftarrow  
\bigoplus\limits_{(a_1, a_2, \cdots, a_{n-1})\in T_{n-1}(E,I) } 
\Lan^S F
\stackrel{d_n}\longleftarrow 
\bigoplus\limits_{(a_1, a_2, \cdots, a_n)\in T_n(E,I) } 
\Lan^S F \leftarrow \cdots
\end{multline*}
where the differentials are defined by (\ref{diff}).
For $(a_1, \ldots, a_n)\in T_n(E,I)$ and $g\in \Lan^S F$, 
the values  $d_n(a_1, \cdots, a_n, g)$ are equal to 
$$
\sum\limits_{s=1}^n (-1)^s (a_1, \cdots, \widehat{a_s} , \cdots, a_n,
\Lan^S F (a_s)(g))
- \sum\limits_{s=1}^n (-1)^s (a_1, \cdots, \widehat{a_s} , \cdots, a_n, g)
$$
Substituting $g=(x,f)\in \bigoplus\limits_{x\in X}F(x)= \Lan^S F$ 
and taking into account the equality 
$\Lan^S F (a_s) (x,f)= (x\cdot a_s, F(x\stackrel{a_s}\rightarrow x\cdot a_s)(f))$ 
realized by Lemma \ref{misc2}, we obtain
\begin{multline*}
d_n(x,a_1, \cdots, a_n,f) = 
\sum_{s=1}^n(-1)^s 
(x\cdot a_s, a_1, \cdots, \widehat{a_s}, \cdots, a_n, 
F(x\stackrel{a_s}\rightarrow x\cdot a_s)(f))
\\ 
-\sum_{s=1}^n(-1)^s 
(x, a_1, \cdots, \widehat{a_s}, \cdots, a_n, f) 
\end{multline*}
\hfill $\eproof$

\begin{corollary}
Suppose that $M(E,I)$ is locally finite-dimensional. Then for any right $M(E,I)$-set $X$,
 the groups  
$\coLim_n^{(M(E,I)/X)^{op}}\Delta\ZZ$ are isomorphic to integer homology groups 
of the semicubical set $Q_* X$.
\end{corollary}

\begin{conjecture}
If $M(E,I)$ is locally finite-dimensional, then for each right $M(E,I)$-set $X$,
 the homology group 
$\coLim_1^{(M(E,I)/X)^{op}}\Delta\ZZ$ is torsion-free.
\end{conjecture}

\begin{example}\label{exstar}
Suppose that $P=\{\star\}$ is the right $M(E,I)$-set over  
a locally finite-dimensional free partially commutative monoid $M(E,I)$. 
By Theorem \ref{mainres}, the groups 
$\coLim_n^{(M(E,I)/P)^{op}}\Delta\ZZ$ are isomorphic 
to the homology groups of the complex
\begin{multline*}
0 \leftarrow 
 \ZZ \stackrel{d_1}\longleftarrow 
 \bigoplus\limits_{(\star, a_1)\in Q_1 P}\ZZ \stackrel{d_2}\longleftarrow 
 \bigoplus\limits_{ 
(\star, a_1,a_2)\in Q_2 P } \ZZ
\leftarrow  \cdots 
\\
\cdots ~~~\leftarrow  
\bigoplus\limits_{(\star, a_1, a_2, \cdots, a_{n-1})\in Q_{n-1}P}\ZZ
\stackrel{d_n}\longleftarrow 
\bigoplus\limits_{(\star, a_1, a_2, \cdots, a_n)\in Q_n P} 
\ZZ \leftarrow \cdots~,
\end{multline*}
where
$d_n(\star, a_1, \cdots, a_n) = 0$. 
Consequently,  $\coLim_n^{(M(E,I)/P)^{op}}\Delta\ZZ\cong \ZZ^{(p_n)}$
where $p_n$ is the cardinality of the set of the subsets
 $\{a_1, \cdots, a_n\} \subseteq E$ consisting of mutually commuting 
elements.
  Here $p_0= 1$ as the number of empty subsets.
\end{example}

\subsection{Homology of asynchronous transition systems}

Following \cite{X20082}, we will consider an asynchronous transition system 
as a  
{\em nondegenerated space of states} with a distinguished {\em initial point}.
Recall that the category of sets and partially functions may be 
considered as the category of pointed sets.

\smallskip
\noindent
{\bf Partially functions and pointed maps.}
A {\em pointed set} $X$ is a set with a distinguished element denoted by $\star$. 
A map $f: X\rightarrow Y$ between pointed sets 
is {\em pointed} if it satisfy  $f(\star)=\star$.

Let ${\rm Set}_*$ be the category of pointed sets and pointed maps between them.
Denote by $\sqcup$ the disjoint union.
The category of sets and maps ${\rm Set}$ admits the inclusion into ${\rm Set}_*$
which assign to every set $S$ the pointed set $S_*= S\sqcup\{\star\}$ and 
to each map $\sigma: S\rightarrow S'$ the map $\sigma_*: S_*\rightarrow S'_*$ defined 
by $\sigma_*(s)=\sigma(s)$ for $s\in S$ and $\sigma_*(\star)=\star$.

For any partially function
$f: X \rightharpoonup Y$ defined in a subset 
$Dom f\subseteq X$, 
we may define the corresponding 
 pointed map $f_*: X_*\rightarrow Y_*$ by
$$
f_*(x)=\left\{
\begin{array}{ll}
f(x) & \mbox{if} \quad  x\in Dom f,\\
{\star} & \mbox{otherwise}.
\end{array}
\right.
$$
By \cite{bed1988}, this correspondence show that 
the category of sets and partially functions is equivalent to  $\Set_*$.

\smallskip
\noindent
{\bf Category of state spaces.}
\begin{definition}
A {\em state space} $\Sigma= (S,E,I,\Tran)$ consists of sets $S$ and $E$,
a subset $\Tran\subseteq S\times E\times S$,
and an irreflexive antisymmetric relation $I\subseteq E\times E$
satisfying to the following two axioms.
\begin{enumerate}
\item if $(s,e,s')\in \Tran$ and $(s,e,s'')\in \Tran$, then $s'=s''$;
\item for any pair $(e_1,e_2)\in I$ and triples $(s,e_1,s_1)\in \Tran$ and
$(s_1,e_2,v)\in \Tran$, there exists  $s_2\in S$ such that
$(s,e_2,s_2)\in \Tran$ and $(s_2,e_1,v)\in \Tran$.
\end{enumerate}
Elements  $s\in S$ are called {\em states},  $(s,e,s')\in \Tran$ 
{\em transitions},  $e\in E$ {\em events}.  $I$ is the
{\em independence relation}.

A state space is called {\em nondegenerate} if it satisfies the condition that 
for every $e\in E$ there exist $s,s'\in S$ such those
$(s,e,s')\in \Tran$.
\end{definition}

Let ${\mathcal S}$ denote the category of state spaces and 
morphisms 
$$
(\sigma,\eta): (S, E, \Tran, I)\rightarrow (S', E', \Tran', I')
$$ 
given by 
maps $\sigma: S\rightarrow S'$ and pointed maps
$\eta: E_*\rightarrow E'_*$ 
satisfying the conditions
\begin{enumerate}
\item for any  
$(s_1, e, s_2)\in \Tran$, 
if $ \eta(e)\not= \star$,  then
$(\sigma(s_1), \eta(e), \sigma(s_2))\in \Tran'$,  otherwise  
$\sigma(s_1)= \sigma(s_2)$;
\item if $(e_1,e_2)\in I$, $\eta(e_1)\not=\star$, 
 $\eta(e_2)\not=\star$,  then  $(\eta(e_1),\eta(e_2))\in I'$.
\end{enumerate}
This definition with the condition $\sigma(s_0)=s'_0$ gives the definition  
of morphisms between asynchronous transition systems.
For $\eta$ satisfying the condition (ii), let 
$\eta_\bullet: E\rightarrow E'\cup \{1\}$ denote the map
$$
\eta_\bullet=\left\{
\begin{array}{ll}
\eta(e), & \mbox{if}   ~  \eta(e)~ \mbox{is defined} , \\
1, & \mbox{otherwise}
\end{array}
\right.
$$
Let $\widetilde{\eta_\bullet}: M(E,I)\rightarrow M(E',I')$ be 
the extension of $\eta_\bullet$ to the monoid homomorphism.

By \cite[Proposition 2]{X20082}, every object $\Sigma$ of  
$\mathcal S$ may be given as a pointed set 
$S_*=S\sqcup\{\star\}$ with a right action of a free partially 
commutative monoid. 
There is a correspondence assigning to any
 morphism from $\Sigma: M(E,I)^{op} \rightarrow \Set_*$ to  
$\Sigma':  M(E',I')^{op} \rightarrow \Set_*$
the diagram

$$
\xymatrix{
	M(E,I)^{\rm op}\ar[rr]^{\widetilde{\eta_\bullet}^{\rm op}} 
 	  \ar[rdd]_{\Sigma}
		&& M(E',I')^{\rm op} 
	\ar[ldd]^{{\Sigma'}}\\
		&{\nearrow}_{\sigma_*}\\
		& {\rm Set}_*
	}
$$
in which  $\sigma_*: \Sigma\rightarrow \Sigma'\circ \widetilde{\eta_\bullet}^{op}$ 
is the natural transformation mapping  
$S$ into $S'$ by $\sigma$. 
It easy to see that the morphism is given by a pair $(\sigma, \varphi)$ 
consisting of a map
$\sigma: S\rightarrow S'$ and a monoid homomorphism
$\varphi: M(E,I)\rightarrow M(E',I')$ such those 
$\varphi(E)\subseteq E'\cup\{1\}$ è $\sigma_*(se)=\sigma_*(s)\varphi(e)$
for all $s\in S$ and $e\in E$.

\smallskip
\noindent
{\bf Homology of state spaces.} Let $U: \Set_*\rightarrow \Set$ be a functor 
which simply forgets the distinguished points. For any state space 
$\Sigma: M(E,I)^{op}\rightarrow \Set_*$ the composition  
$U\circ\Sigma$ is a right $M(E,I)$-set. 
Denote by $K_*(\Sigma)$ the category $(M(E,I)/(U\circ\Sigma))^{op}$.
Its objects may be considered as elements in  $S_*$, and morphisms are 
triples $s_1\stackrel{\mu}\rightarrow s_2$, $\mu\in M(E,I)$, $s_1\in S_*$, 
$s_2\in S_*$. Composition of morphisms  
$(s_2\stackrel{\mu_2}\rightarrow s_3)\circ(s_1\stackrel{\mu_1}\rightarrow s_2)$
equals $(s_1\stackrel{\mu_1\mu_2}\rightarrow s_3)$.

{\em Homology groups of a state space with 
coefficients in a functor
$F: K_*(\Sigma)\rightarrow \Ab$} are defined by 
$H_n(\Sigma, F)=\coLim_n^{K_*(\Sigma)}F$ for $n\geq 0$.

{\em Semi-regular higher-dimensional automata} \cite{gou1995} 
are precisely the semicubical sets. So for such automata $X$, there are defined homology 
groups $H_n(X,F)$ with coefficients in homological systems 
$F$. It follows from Theorem \ref{mainres} the following assertion.
\begin{corollary}\label{corres}
Let $\Sigma=(S, E, I, \Tran)$ be a state space.
If
$M(E,I)$ is locally finite-dimensional, then 
$H_n(\Sigma, F)\cong H_n(Q_*(U\circ\Sigma), \overline{F})$ for all $n\geq 0$.
\end{corollary}

\noindent
{\bf Direct summands of integer homology groups.}
We will show that the homology groups of the one-point set considered in Example \ref{exstar} 
are direct summands of  
$H_n(\Sigma,\Delta\ZZ)$.
To that end we define functors with values  $\ZZ$ 
and $0$ on objects.
For a small category $\mC$, a subset $S\subseteq \Ob\mC$ is called
 {\em closed in $\mC$} if it contains with any 
$s\in S$ all objects $c\in \Ob\mC$ those admit morphisms $s\rightarrow c$.
For example, the subset $\{\star\}$ is closed in  $K_*(\Sigma)$. 
The complement of a closed subset in $\Ob\mC$ is called {\em open}. If a set 
$S\subseteq \Ob\mC$ is equal to the intersection of an open subset  and a closed subset, 
then we can define a functor  
$\ZZ[S]$ by
$$
\ZZ[S](c)=\left\{
\begin{array}{ll}
\ZZ & \mbox{if}   ~~  c\in S , \\
0, & \mbox{otherwise,}
\end{array}
\right.
$$
on objects $c\in \Ob\mC$. We put $\ZZ[S](c_1\rightarrow c_2)=1_{\ZZ}$ if $c_1\in S$ and $c_2\in S$,
and $\ZZ[S](c_1\rightarrow c_2)=0$ on the other morphisms.

\begin{proposition}\label{sphom}
Suppose that $\Sigma=(S,E,I,\Tran)$ is an arbitrary state space. 
Then 
$H_n(\Sigma,\Delta\ZZ)\cong \ZZ^{(p_n)}\oplus H_n(\Sigma, \ZZ[S])$ for all $n\geq 0$.
\end{proposition}
{\sc Proof.} 
Consider the full subcategory $K_*(\emptyset)\subseteq K_*(\Sigma)$ 
with  $\Ob(K_*(\emptyset))= \{\star\}$. The inclusion of this full subcategory is 
a coretraction. 
The exact sequence 
$0 \rightarrow \ZZ[\star]\rightarrow \Delta\ZZ \rightarrow \ZZ[S]\rightarrow 0$ 
in $\Ab^{K_*(\Sigma)}$ gives the exact sequence of complexes
\begin{equation}\label{exactsp}
0\rightarrow C_*(K_*(\Sigma), \ZZ[\star])\rightarrow C_*(K_*(\Sigma), \Delta\ZZ)
\rightarrow C_*(K_*(\Sigma), \ZZ[S])\rightarrow 0
\end{equation}
The chain homomorphism  
$C_*(K_*(\Sigma), \ZZ[\star])\rightarrow C_*(K_*(\Sigma), \Delta\ZZ)$ 
is equal to the composition of the isomorphism   
$C_*(K_*(\Sigma), \ZZ[\star])\rightarrow C_*(K_*({\emptyset}), \Delta\ZZ)$ 
and coretraction
$C_*(K_*({\emptyset}), \Delta\ZZ)\rightarrow C_*(K_*(\Sigma), \Delta\ZZ)$.
Hence the exact sequence  (\ref{exactsp})
splits.
The corresponding exact sequence of $n$-th homology groups gives the required
assertion.
\hfill $\Box$

In particular,
if $H_n(\Sigma, \Delta\ZZ)\cong H_n(\pt)$ for all $n\geq 0$, 
then $p_n=0$ for $n>0$.
In this case $E=\emptyset$, $I=\emptyset$, and hence $K_*(\Sigma)$ is a discrete category.
 Since $H_0(K_*(\Sigma))=\ZZ$, this category has an unique object.
Consequently, $S=\emptyset$. It follows the following assertion important for 
the calassification of state spaces.
\begin{corollary}
Let $\Sigma=(S,E,I, \Tran)$ be a state space.
If $H_n(\Sigma, \Delta\ZZ)=0$ for all $n>0$ and $H_0(\Sigma, \Delta\ZZ)=\ZZ$, 
then $S=E=I= \Tran= \emptyset$.
\end{corollary}

\noindent
{\bf Computing homology groups the state space 
consisting of an unique element.}
It may be seemed that the homology groups of state spaces are torsion-free.  
We will show that this opinion is wrong.

Recall that a
{\em simplicial scheme} is a pair $(X,\mathfrak{M})$ consisting of a set 
 $X$ and a set $\mathfrak{M}$ of its finite nonempty subsets satisfying the 
following conditions
\begin{enumerate}
\item $x\in X \Rightarrow \{x\}\in \mathfrak{M}$
\item $S \subseteq T \in \mathfrak{M} \Rightarrow S \in \mathfrak{M}$.
\end{enumerate}
Let $(X,\mathfrak{M})$ be a simplicial set.
Consider an arbitrary total ordering $<$ on $X$.
Denote 
$$
X_n = \{ (x_0, x_1, \cdots, x_n) : \{x_0, x_1, \cdots, x_n\} 
\in \mathfrak{M}
\quad \& \quad x_0< x_1< \cdots< x_n \}.
$$
For any set $S$, let $L(S)$ be a free Abelian group generated by $S$. 
Consider a family of Abelian groups
$$
C_n(X,\mathfrak{M})=\left\{
\begin{array}{ll}
L(X_n) & $for$ \quad  n \geq 0,\\
0, & $otherwise$.
\end{array}
\right.
$$
Define homomorphisms 
$d_n: C_n(X,\mathfrak{M}) \rightarrow C_{n-1}(X,\mathfrak{M}) $
on the generators by 
$$
d_n(x_0, x_1, \cdots, x_n) = 
\sum_{i=0}^n (-1)^i (x_0, \cdots, x_{i-1}, x_{i+1}, \cdots, x_n).
$$

It is well known that the family $(C_n(X,\mathfrak{M}), d_n)$ 
is a chain complex.
Denote  $C_*(X,\mathfrak{M})=(C_n(X,\mathfrak{M}), d_n)$.
We can prove 
that homology groups of the complex $C_*(X,\mathfrak{M})$ 
does not depend on the total ordering of $X$. 
They are called the {\em homology groups  
$H_n(X,\mathfrak{M})$
of the simplicial scheme $(X,\mathfrak{M})$}.

\begin{theorem}\label{exsimpsh}
Let $\Sigma= (\{x_0\}, E, I, \Tran)$ be a state space.
Suppose that the action of $M(E,I)$ is defined by 
$x_0\cdot a= \star$ for every $a\in E$.
Let 
$(E, \mathfrak{M})$ be a simplicial scheme where $\mathfrak{M}$ 
consists of nonempty finite subsets of mutually 
commuting generators and let $p_n=\vert E_{n-1}\vert$ 
denote the cardinal number of $n$-cliques 
in the independence graph of $M(E,I)$. 
If $M(E,I)$ is locally finite-dimensional, then  
 $H_n(\Sigma, \Delta\ZZ) \cong \ZZ^{(p_n)}\oplus H_{n-1}(E,\mathfrak{M})$
for all $n\geq 2$. 
\end{theorem}
{\sc Proof.}
 Let 
$\ZZ[x_0]: K_*(\Sigma) \rightarrow \Ab$ be the functor with values 
 $\ZZ[x_0](x_0)=\ZZ$ and $\ZZ[x_0](\star)=0$.
By Theorem \ref{mainres}, the homology groups 
$H_n(\Sigma, \ZZ[x_0])= \coLim_n^{(M(E,I)/(U\circ\Sigma))^{op}}\ZZ[x_0]$ 
are isomorphic to the homology groups of the complex
\begin{multline*}
0 \leftarrow 
 \ZZ \stackrel{d_1}\longleftarrow 
 \bigoplus\limits_{a_1\in E_1}\ZZ \stackrel{d_2}\longleftarrow 
 \bigoplus\limits_{(a_1,a_2)\in E_2 } \ZZ
\leftarrow  \cdots 
\\
\cdots \leftarrow  
\bigoplus\limits_{(a_1, a_2, \cdots, a_{n-1})\in E_{n-1} } 
\ZZ
\stackrel{d_n}\longleftarrow 
\bigoplus\limits_{(a_1, a_2, \cdots, a_{n})\in E_{n} } 
\ZZ \leftarrow \cdots~,
\end{multline*}
with the differentials
$d_n(a_1,\cdots, a_n) = 
\sum_{s=1}^n(-1)^{s+1} 
(a_1, \cdots, \widehat{a_s}, \cdots, a_n)$.
This complex equals the shifted complex $(C_{k-1}(E, \mathfrak{M}), d_{k-1})$
 in the dimensions $k\geq 1$. Hence its  $k$-th homology groups 
are isomorphic to $H_{k-1}(E,\mathfrak{M})$ for  $k\geq 2$.

By Proposition \ref{sphom}, we obtain
$H_n(\Sigma, \Delta\ZZ)\cong \ZZ^{(p_n)}\oplus H_n(\Sigma, \ZZ[x_0])$ for all 
$n\geq 0$.
The required assertion follows from $H_n(\Sigma, \ZZ[x_0])\cong H_{n-1}(E,\mathfrak{M})$ for 
all $n\geq 2$.  
\hfill $\eproof$

Any asynchronous transition system may be considered 
as a pair $T=(\Sigma,s_0)$ of its nondegenerate state space 
$\Sigma$ with an initial state $s_0\in S$.
Each morphism between asynchronous transition systems 
$(\Sigma, s_0)\rightarrow (\Sigma', s_0')$ may be given by a
morphism of state spaces  $(\sigma, \eta): \Sigma\rightarrow \Sigma'$ such that
 $\sigma(s_0)=s_0'$.

Let $T=(\Sigma, s_0)$ be an asynchronous transition system with a state space  
$\Sigma=(S, E, I, \Tran)$.
Events $e_j\in E$, $j\in J$, are caled {\em mutually independent} if
 $(e_j,e_{j'})\in I$ for all $j,j'\in J$ such those $j\not=j'$.
A state $s\in S\sqcup\{\star\}$ is {\em available} if there exists  
$\mu \in M(E,I)$ such that
$s_0\cdot \mu = s$. 
The set of available states is denoted by $S(s_0)$.
The monoid $M(E,I)$ will act on the set of the available states.
Suppose that $T(s_0)$ is the state space whose set of states equals $S(s_0)$.

Define {\em homology groups of space of available states} 
by $H_n(K_*(T(s_0)), \Delta\ZZ)$. 
Applying Theorem \ref{mainres}, we get the following.

\begin{corollary}\label{main3}
Suppose that an asynchronous transition system $T=(\Sigma, s_0)$ 
 does not contain infinite subsets of mutually 
independent events. Let 
$\Sigma=(S, E, I,\Tran)$ be its the state space. 
Suppose that the set $E$ is totally ordered.
Then $H_n(K_*(T(s_0)), \Delta\ZZ)$ are isomorphic to 
homology groups of the complex 
\begin{multline}\label{complex}
0 \leftarrow 
\bigoplus\limits_{s\in S(s_0)} \ZZ \stackrel{d_1}\longleftarrow 
\bigoplus\limits_{(s,e_1)\in S\times E_1} \ZZ
\stackrel{d_2}\longleftarrow \bigoplus\limits_{{(s, e_1, e_2)\in S\times E_2}} \ZZ
\leftarrow  \cdots 
\\
\cdots ~~
\leftarrow  
\bigoplus\limits_{(s, e_1,\cdots, e_{n-1})\in S\times E_{n-1}}
\ZZ
\stackrel{d_n}\longleftarrow 
\bigoplus\limits_{(s, e_1,\cdots, e_{n})\in S\times E_{n}} 
\ZZ \leftarrow \cdots~ 
\end{multline}
where $E_n$ consists of tuples $e_1<\cdots<e_n$ of mutually commuting elements of
 $E$ with
$$
d_n(s,e_1\cdots e_n) = 
\sum_{i=1}^n(-1)^i \left(
(s\cdot e_i, e_1, \cdots, \widehat{e_i}, \cdots, e_n)- 
(s, e_1, \cdots, \widehat{e_i}, \cdots, e_n) \right)
$$ 
\end{corollary}
If the cardinal numbers of mutually independent events is bounded above by a natural number,
then this comlex has a finite length.

\section{Concluding remarks}

Let $\Sigma=(S,E,I,Tran)$ be a state space.
Consider the full subcategory $K(\Sigma)\subset K_*(\Sigma)$ whose objects are 
elements $s\in S$. In \cite{X20042}, it was studied homology groups of the category 
 $K(\Sigma)$. It was built an algorithm of computing the first integer homology group of 
this category and applied for the calculation of homology groups of finite Petri CE nets.
An algorithm of computing all integer homology groups 
of finite Petri CE nets is not found \cite[Open Problem 1]{X20042}.
We put forward a conjecture whose confirmation would be to solve this problem.
Let $Q_*'(U\circ\Sigma)\subseteq Q_*(U\circ\Sigma)$ be 
a semicubical subset
 consisting of sets 
$$
Q_n'(U\circ\Sigma)= \{(s, e_1,\cdots, e_{n})\in Q_n(U\circ\Sigma): 
s e_1 \cdots e_{n}\not=\star\}
$$
Consider any functor
 $F: K(\Sigma)\rightarrow Ab$. Extend it to $K_*(\Sigma)$ by
 $F(\star)=0$.

\begin{conjecture}
Let $\Sigma=(S,E,I,Tran)$ be a state space. If the monoid
 $M(E,I)$ is locally finite-dimensional, then for all integer $n\geq 0$,
 the groups
$\coLim_n^{K(\Sigma)}F$ are isomorphic to  $n$-th homology groups of 
the complex consisting of groups 
$\bigoplus\limits_{(s, e_1,\cdots, e_{n})\in Q'_n(U\circ\Sigma) }
F(s)$ and differentials given by 
$d_n(s,e_1, \cdots, e_n,f) =$
\begin{multline*}
\sum_{i=1}^n(-1)^i (
(s\cdot e_i, e_1, \cdots, \widehat{e_i}, \cdots, e_n, 
F(s\stackrel{e_i}\rightarrow s\cdot a_i)(f))\\
- (s, e_1, \cdots, \widehat{e_i}, \cdots, e_n, f) )
\end{multline*}

\end{conjecture}


\begin{thebibliography}{1}

\bibitem{nie1996}
Nielsen M., Winskel G., 
        \newblock{``Petri nets and bisimulation''},
        \newblock{ {\it Theoretical Computer Science}, {\bf 153}:1-2, 
        1996, 211--244}
\bibitem{shi1985}
Shields M.W.,
        \newblock{``Concurrent machines'',}
        \newblock{ {\it Computer Journal}, {\bf 28} (1985), 449--465}
\bibitem{bed1988}
Bednarczyk M. A., 
	\newblock{\it Categories of Asynchronous Systems}, 
	\newblock{Ph.D. thesis, University of Sussex, report 1/88, 1988;
		http://www.ipipan.gda.pl/$~\widetilde{~}$marek
       }
\bibitem{X20032}
{Husainov A. A., Tkachenko V. V.,}
     \newblock{  ``Homology groups of asynchronous transition
systems'',}
     \newblock{ {\it Mathematical modeling and the near questions of mathematics.
Collection of the scientifics works.}, KhGPU, Khabarovsk, 2003,
     23--33 (Russian)
}

\bibitem{X20042}
{Husainov A.},
        \newblock{``On the homology of small categories and asynchronous 
	transition systems''},
        \newblock{
        {\it Homology Homotopy Appl.}, {\bf 6}:1 (2004), 439--471;
       http://www.rmi.acnet.ge/hha
        }
\bibitem{gou1995}
{Goubault E.},
     \newblock{\it  The Geometry of Concurrency},
     \newblock{Ph.D. Thesis, Ecole Normale Sup\'erieure, 1995;
     http://www.dmi.ens.fr/$~\widetilde{~}$goubault}
\bibitem{gau2000}
{Gaucher P.},
     \newblock{``About the globular homology of higher
     dimensional automata'',}
     \newblock{ {\it Cah. Topol. Geom. Differ.}, {\bf 43}:2, 2002, 107--156}
\bibitem{X20082}
Khusainov A.A., Lopatkin V.E., Treshchev I.A.,
     \newblock{``Algebraic topology approach to mathematical model analysis 
	of concurrent computational processes'',}
     \newblock{ {\it Sibirski\u{\i} Zhurnal Industrial'no\u{\i} Matematiki}, {\bf 11}:1 (2008), 
	 141--152 (Russian)}
\bibitem{pol2007}
Polyakova L.Yu.,
     \newblock{``Resolutions for free 
	partially commutative monoids'',}
     \newblock{ {\it Sib. Math. J.}, {\bf 48}:6 (2007), 1038--1045; 
	translation from {\it Sib. Mat. Zh.}, {\bf 48}:6 (2007), 1295--1304}
\bibitem{X20081}
Khusainov A.A.,
     \newblock{``Homology groups of semicubical sets'',}
     \newblock{{\it Sib. Math. J.}, {\bf 49}:1 (2008), 180--190;
	translation from {\it Sib. Mat. Zh.}, {\bf 49}:1 (2008), 224--237}
\bibitem{kac2004}
Kaczynski T., Mischaikov K., Mrozek M.,
        \newblock{\it Computational homology},
        \newblock{ Springer-Verlag, New York, 2004, (Appl. Math. Sci.; 157)}
\bibitem{gab1967}
{Gabriel P., Zisman M.,}
     \newblock{\it Calculus of Fractions and Homotopy Theory }
     \newblock{ [Russian translation], Moscow, Mir, 1971}
\bibitem{mac1998}
MacLane S.,
        \newblock{\it Categories for the Working Mathematician,}
        \newblock{Springer-Verlag, Berlin, 1971}
\bibitem{obe1968}
Oberst U.,
        \newblock{``Homology of categories and exactness
        of direct limits'',}
        \newblock{ {\it Math. Z.}, {\bf 107} (1968), 87--115}
\bibitem{X20083}
Husainov A. A., 
        \newblock{``On the Leech dimension of a free
        partially commutative monoid''},
        \newblock{ {\it Tbilisi Math. J.}, {\bf 1}:1 (2008), 71--87;
	http://ncst.org.ge/Journals/TMJ/index.html
}
\bibitem{kac2003}
Kaczynski T., Mischaikov K., Mrozek M.,
     \newblock{``Computational homology''},
     \newblock{ {\it Homology Homotopy Appl.}, {\bf 5}:2 (2003), 233 -- 256}
\bibitem{bau1985}
Baues H.-J., Wirsching G.,
        \newblock{``Cohomology of small categories'',}
        \newblock{ {\it J. Pure Appl. Algebra}, {\bf 38}:2-3 (1985), 187--211}
\bibitem{die1997}
{Diekert V., M{\'e}tivier Y.,}
     \newblock{``Partial Commutation and Traces'',}
     \newblock{{\it Handbook of formal languages}, {\bf 3}, Springer-Verlag, 
	New York, 1997, 457--533}


\end{thebibliography}
\end{document}